\def\rest{\mathord{\restriction}}

\def\phi{\varphi}

\def\su{\subseteq}
\def\a{\alpha}
\def\b{\beta}
\def\g{\gamma}

\def\d{\delta}
\def\l{\lambda}

\def\z{\zeta}

\def\om{\omega}

\def\lng{\langle}
\def\rng{\rangle}
\def\ov{\overline}

\def\cl{{\text {cl}\,}}
\def\dom{{\text dom\,}}
\def\cf{{\text {cf}\,}}

\def\pcf{{\text {pcf}}}

\def\imply{\Rightarrow}

\outer\long\def\ignore#1\endignore{}

\newcount\itemno
\def\itm{\advance\itemno1 \item{(\number\itemno)}}

\def\aitm{\advance\itemno1
\item{(\letter\itemno)}}

\def\letter#1{\ifcase#1 \or a\or b\or c\or d\or e\or f\or g\or h\or
i\or j\or k\or l\or m\or n\or o\or p\or q\or r\or s\or t\or u\or v\or
w\or x\or y\or z\else\toomanyconditions\fi}
\def\raitm{\advance\itemno1 \item{(\rletter\itemno)}}
\def\rletter#1{\ifcase#1\or `\or a\or b\or c\or d\or e\or f\or g\or
h\or i\or k\or k\or l\or n\or n\or q\or r\or t\or v
\else\toomanyconditions\fi}

\documentclass[12pt]{article}
\usepackage{amstext}
\usepackage{amssymb}
\usepackage{amsfonts}
\usepackage{amsthm}
\usepackage{amsopn}

\newtheorem{theorem}{Theorem}
\newtheorem{lemma}{Lemma}

 \newtheorem{definition}[lemma]{Definition}

 \newtheorem{problem}[lemma]{Problem}
 \newtheorem{claim}[lemma]{Claim}

\begin{document}
\author{Menachem Kojman\\
Department Of Mathematics\\
Carnegie-Mellon University
\and Saharon Shelah\thanks{Research supported by  ``The Israel
Science Foundation'' 	administered by The Israel Academy of Sciences
and Humanities. Publication 609.  }\\
Institute of Mathematics\\
The Hebrew University of Jerusalem\\
and \\
Department of Mathematics\\
Rutgers University}
\title{A ZFC Dowker space in $\aleph_{\om+1}$: an
application of pcf theory to topology}
\date{October 95}
\maketitle

\begin{abstract} The existence of a Dowker space of  cardinality
$\aleph_{\om+1}$ and weight $\aleph_{\omega+1}$ is proved in ZFC using
pcf theory.
\end{abstract}

\section{Introduction}
A Dowker space is a normal Hausdorff topological space whose product
with the unit interval is not normal. The problem of existence of such
spaces was raised by C.~H.~ Dowker in 1951. C.~H.~Dowker characterized Dowker
spaces as normal Hausdorff and not countably paracompact \cite{dowker}.

Exactly two Dowker spaces were constructed in ZFC so far.  The existence
of a Dowker space in ZFC was first proved by M.~E.~Rudin in 1971
\cite{rudin}, and her space was the only known Dowker space in ZFC for
over two decades. Rudin's space is a subspace of $\prod_{n\ge
1}(\aleph_n+1)$ and has cardinality $\aleph_\om^{\aleph_0}$. The
problem of finding a Dowker space of smaller cardinality in ZFC was
referred to as the ``small dowker space problem''.

Z.~T.~Balogh constructed recently \cite{balogh} a dowker space in ZFC
whose cardinality is $2^{\aleph_0}$.

While both Rudin's and Balogh's spaces are constructed in ZFC, their
respective cardinalities are not decided in ZFC, as is well known by
the independence results of P.~Cohen: both $2^{\aleph_0}$ and
$\aleph_\om^{\aleph_0}$ have no bound in ZFC, (and may be equal to
each other).

The problem of which is the first $\aleph_\a$ in which ZFC proves the
existence of a Dowker space remains thus unanswered by Rudin's and
Balogh's results.

In this paper we prove that there is a Dowker space of
cardinality $\aleph_{\om+1}$. A non-exponential bound is thus provided
for the cardinality of the smallest ZFC Dowker space. We do this by
exhibiting a  Dowker subspace of Rudin's space of that cardinality. Our
construction avoids the exponent which appears in the cardinality of
Rudin's space by working with only a fraction of
$\aleph_\om^{\aleph_0}$. 
It remains open whether $\aleph_{\om+1}$ is the {\em first} cardinal at
which there is a ZFC Dowker space. 

We shall describe shortly the cardinal arithmetic developments which
enable this result. The next three paragraphs are not necessary for
understanding the proofs in this paper.

In the last decade there has been a considerable
advance in understanding of the infinite exponents of singular
cardinals, in particular the exponent $\aleph_\om^{\aleph_0}$. This
exponent is the product of two factors: $2^{\aleph_0}\times
\cf\lng[\aleph_\om]^{\aleph_0},\su\rng$. The second factor, the cofinality
of the partial ordering of inclusion over all countable subsets of
$\aleph_\om$, is the least number of countable subsets of $\aleph_\om$
needed to cover every countable subset of $\aleph_\om$; the first factor
is the number of subsets of a single countable set. Since
$\aleph_\om^{\aleph_0}$ is the number of countable subsets of
$\aleph_\om$, the equality $\aleph_\om=2^{\aleph_0}\times \cf\lng
[\aleph_\om]^{\aleph_0}\rng$ is clear.

While for $2^{\aleph_0}$ it is consistent with ZFC to equal any
cardinal of uncountable cofinality, the second author's work on
Cardinal Arithmetic provides  a ZFC bound of $\aleph_{\om_4}$ on  the
factor $\cf\lng [\aleph_\om]^{\aleph_0},\su\rng$. 

This is done  by approximating $\cf\lng
[\aleph_\om]^{\aleph_0},\su\rng$ by an interval of regular cardinals,
whose first element is $\aleph_{\om+1}$ and whose last element is
$\cf\lng [\aleph_\om]^{\aleph_0},\su\rng$, and so that every regular
cardinal $\l$ in this interval is the {\em true cofinality} of a
reduced product $\prod B_\l/J_{<\l}$ of a set $B_\l\su
\{\aleph_n:n<\om\}$ modulo an ideal $J_{<\l}$ over $\om$. 
The theory of reduced products of small sets of regular cardinals,
known now as {\em pcf theory}\footnote{$\pcf$ means {\em possible
cofinalities}}, is used to put a bound of $\om_4$ on the length of
this interval.

Back to topology now, it turns out that the pcf approximations to
$\aleph_{\om}^{\aleph_0}$ are concrete enough to ``commute'' with
Rudin's construction of a Dowker space. Rudin defines a topology on a
subspace of the functions space $\prod_{n>1}(\aleph_n+1)$. What is
gotten by restricting Rudin's definition to the first approximation of
$\aleph_{\om}^{\aleph_0}$ is a {\em closed} and {\em cofinal} Dowker
subspace $X$ of the Rudin space $X^R$ of cardinality $\aleph_{\om+1}$.

Hardly any background is needed to state the pcf theorem we are using
here.  However, an interested reader can find presentations of pcf
theory in either \cite{burk:mag}, the second author's \cite{thebook}
or the first author's \cite{pcfPap}. The pcf theorem used here is
covered in detail in each of those three sources.

\medskip
\noindent
\textbf{Acknowledgements} 
The first author wishes to thank J.~Baumgartner, for inviting him to
the 11th Summer Conference on General Topology and Applications in
Maine in the summer of 95, Z.~T.~Balogh, for both presenting the small
Dowker space problem in that meeting and for very pleasant and
illuminating conversations that followed, and P.~Szeptycki, for
suggesting that the subspaces we construct may actually be closed, and 
Mirna Djamonja, for  detecting  an inaccuracy in an earlier version of
this paper.

\section{Notation and pcf}

In this section we present a few simple definitions needed to state
the pcf theorem used in proving the
existence of an $\aleph_{\om+1}$-Dowker space.

Suppose $B\su \om$ is a subset of the natural numbers.

\begin{definition}

\begin{enumerate}
\item $\prod_{n\in B}\aleph_n=\{f: \dom f=B\wedge f(n)<\aleph_n \text{ for }
n\in B\}$ 
\item $\prod_{n\in B}(\aleph_n+1)=\{f: \dom f=B\wedge f(n)\le \aleph_n
\text{ for }  
n\in B\}$

\item for  $f,g\in \prod _{n\in B}(\aleph_n+1)$ let:
\begin{enumerate}
\item $f<g$ iff $\forall n\in B\left[f(n)<g(n)\right]$
\item $f\le g$ iff $\forall n\in B\left[f(n)\le g(n)\right]$
\item  $f\le^*g$ iff
$\{n:f(n)>g(n)\}$ is finite
\item  $f<^* g$ iff $\{n:f(n)\ge g(n)\}$
is finite 
\item  $f=^*g$ iff $\{n:f(n)\not= g(n)\}$ is finite
\end{enumerate}

\item A sequence $\lng f_\a:\a<\l\rng$ of functions in $\prod_{n\in
B}\aleph_n$ is {\em increasing in}  $<$ ($\le,<^*,\le^*$) iff
$\a<\b<\l\imply f_\a<f_\b$ ($f_\a\le f_\b$, $f_\a<^* f_\b$,
$f_\a\le^* f_\b$)

\item $g\in \prod _{n\in B}(\aleph_n+1)$ is an {\em upper bound}
of $\{ f_\a:\a<\d\}\su \prod_{n\in B}\aleph_n$ if and only if
$f_\a\le^* g$ for all $\a<\d$

\item $g\in \prod _{n\in B}(\aleph_n+1)$ is a {\em least upper bound}
of $\{ f_\a:\a<\d\}\su \prod_{n\in B}\aleph_n$ if and only if $g$ is an
upper bound of $\{ f_\a:\a<\d\}\su \prod_{n\in B}\aleph_n$ and if $g'$
is an upper bound of $\{ f_\a:\a<\d\}$ then
$g\le^* g$
\end{enumerate}
\end{definition}

\begin{theorem} (Shelah)  \label{tcf} There is a set $B=B_{\aleph_{\om+1}}\su
\om$ and a sequence $\ov f=\lng f_\a:\a<{\aleph_{\om+1}}\rng$ of
functions in  $\prod_{n\in B}\aleph_n$ such that:
\begin{itemize}
\item $\ov f$ is increasing in $<^*$
\item $\ov f$ is cofinal: for every $f\in \prod_{n\in B}\aleph_n$
there is $\a<{\aleph_{\om+1}}$ so that $f<^* f_\a$
\end{itemize}
\end{theorem}

A sequence as in the theorem above will be referred to as an
``$\aleph_{\om+1}$-scale''.

By Theorem \ref{tcf} we can find $B\su \om$ and an
$\aleph_{\om+1}$-scale $\ov g=\lng g_\a:\a<\aleph_{\om+1}\rng$ in
$\prod _{n\in B}\aleph_n$. The set $B$ is clearly infinite. Restricting
every $g_\a\in \ov g$ to a fixed co-finite set of coordinates does
not matter, so we assume without loss of generality that $0,1\notin B$.
For notational simplicity we pretend that $B=\om-\{0,1\}$; if this is not
the case, we need to replace $\aleph_n$ in what follows by the $n$-th
element of $B$. We sum up our assumptions in the following:

\begin{claim} \label{ass} We can assume without loss of generality
that there is an $\aleph_{\om+1}$-scale $\ov g=\lng
g_\a:\a<{\aleph_{\om+1}}\rng$ in $\prod_{n>1}\aleph_n$.
\end{claim}

\begin{claim} \label{cont}  There is an $\aleph_{\om+1}$-scale $\ov
f=\lng f_\a:\a<{\aleph_{\om+1}}\rng$ in
$\prod_{n>1}\aleph_n$ so that for every $\d<{\aleph_{\om+1}}$, if
$\cf\d>\aleph_0$ and a least upper bound of $\ov f\rest \d$ exists,
then $f_\d$ is a least upper bound of $f\rest \d$. 
\end{claim}

\begin{proof} Fix an $\aleph_{\om+1}$-scale $\ov g=\lng
g_\a:\a<{\aleph_{\om+1}}\rng$ in $\prod_{n>1}\aleph_n$ as guranteed by
Claim \ref{ass}.  Define $f_\a$ by induction on $\a<{\aleph_{\om+1}}$
as follows: If $\a$ is successor or limit of countable cofinality let
$f_\a$ be $g_\b$ for the first $\b\in (\a,{\aleph_{\om+1}})$ for which
$g_\b>^* f_\b$ for all $\b<\a$. If $\cf \a>\aleph_0$ then let $g_\a$
be a least upper bound to $\ov f\rest \d:=\lng f_\b:\b<\a\rng$, if
such least upper bound exists; else, define $f_\a$ as in the previous
cases.

The sequence $\ov f=\lng f_\a:\a<{\aleph_{\om+1}}\rng$ is increasing
cofinal in $\prod_{n>1}\aleph_n$ and by its definition satisfies the
required condition.
\end{proof}

\begin{claim}\label{suf} Suppose $0<m\le k<\om$. Let $\lng
\a(\z):\zeta<\aleph_m\rng$ be strictly  increasing with
$\sup\{\a(\z):\z<\aleph_m\}=\d<{\aleph_{\om+1}}$. 
If $\lng g_\zeta:\zeta<\aleph_m\rng$ is a sequence of functions in
$\prod_{n> k}\aleph_n$ which is increasing in $<$, and  $g_\z=^*
f_{\a(\z)}$ for every $\zeta<\aleph_m$, then:
\begin{itemize}
\item $g:=\sup \{g_\z:\z<\aleph_m\}\in \prod _{n>k}\aleph_n$ is a least
upper bound of $\ov f\rest \d$ 
\item $\cf g(n)=\aleph_m$ for all $n> k$
\item $g=^* f_\d$
\end{itemize}
\end{claim}

\begin{proof}[Proof of Claim] 

Let $g:=\sup \{g_\zeta:\z<\aleph_m\}$. Since $\lng g_\z:\z<\aleph_m\rng$
is increasing in $<$, necessarily $\cf g(n)=\aleph_m$ for all $n\ge k$
and since $g(n)\le \aleph_n$ it follows that $g(n)<\aleph_n$ for
$n>k$  and therefore $g\in
\prod _{n>k}\aleph_n$.

Suppose that $\gamma<\d$ is
arbitrary. There exists $\z<\aleph_m$ such that $\gamma<\a(\zeta)$,
hence $f_\gamma<^* f_{\a(\z)}=^*g_\z\le g$. Thus $g$ is an upper bound
or $\ov f\rest \d$. 

To show that $g$ is a least upper bound suppose that $g'$ is an upper
bound of $\ov f\rest \d$. Let $X:=\{n> k: g'(n)<g(n)\}$. For every
$n\in X$ find $\zeta(n)<\aleph_m$ such that
$g_{\zeta(n)}(n)>g'(n)$. Such $\zeta(n)$ can be found because $g=\sup
\{g_\z:\z<\aleph_m\}$. Let $\zeta^*:=\sup\{\zeta(n):n>1\}$. Since
$\aleph_m>\aleph_0$, $\zeta^*<\aleph_m$. Since $\lng
g_\z:\z<\aleph_m\rng$ is increasing in $<$, it holds
that $f_{\z(*)}\ge f_{\z(n)}(n)>g'(n)$ for every $n\in X$ . But $g'$ is
an upper bound of $\ov f\rest\d$, so
$f_{\z(*)}\le^* g'$ and $X$ is therefore finite. 

By the definition of $\ov f$ we conclude that $f_\d$ is a least upper
bound of $\ov f\rest \d$. Since both $g$ and $f_\d$ are least upper
bounds of $\ov f\rest \d$ it follows that $g=^* f_\d$.
\end{proof}

\section{The Space}

\begin{definition}\label{xrudin}  Let $X^R=\{h\in
\prod_{n>1}(\aleph_n+1): \exists m\,\forall 
n\,\left[\aleph_0<\cf h(n)<\aleph_m\right]\}$.
\end{definition}

The space is $X^R$ is {\em the Rudin space} from \cite{rudin}
 with the Hausdorff topology defined by letting, for every $f<g$ in
$\prod_{n>1}(\aleph_n+1)$, 
\begin{equation}\label{top} (f,g]:=\{h\in X^R:f<h\le g\}
\end{equation}
 
be a basic open set (see \cite{rudin}). 

Recall that a normal Hausdorff space is {\em countably
paracompact} iff for every decreasing sequence $\lng D_n:n<\om\rng$ of
closed sets such that $\bigcap D_n=\emptyset$ there are open sets
$U_n\supseteq D_n$ with $\bigcap U_n=\emptyset$.

\begin{definition} \label{dn} $D_n:=\{h\in X^R:\exists m\,\ge n \,\bigl[
h(m)=\aleph_m\bigr]\}$
\end{definition}

M.~E.~Rudin defined in \cite{rudin} the closed subsets $D_n\su X^R$
above and proved:

\begin{theorem}(Rudin)\label{xrdowker}
\begin{enumerate} 
\item $X^R$ is collectionwise normal
\item If $U_n\su X^R$ is open and $D_n\su U_n$ for all $n>1$ then
$\bigcap_{n>1} U_n$ is not empty
\end{enumerate}
\end{theorem}

Those two facts establish by \cite{dowker} that $X^R$  Dowker.

Let $\ov f=\lng f_\a:\a<{\aleph_{\om+1}}\rng$ be as provided by Claim
\ref{cont}. We use this scale to extract a closed Dowker subspace of
cardinality ${\aleph_{\om+1}}$ from Rudin's space.

\begin{definition}\label{X}
 $X=\{h\in X^R:\exists \a<\aleph_{\om+1}\,[h=^* f_\a]\}$
\end{definition}

Since $|\{h\in X^R:h=^* f_\a\}|=\aleph_\om$ for every
$\a<\aleph_{\om+1}$  it is obvious that
$|X|={\aleph_{\om+1}}$.

Since 
$\ov f$ is totally ordered by $<^*$, for every $h\in X$ there exists a
{\em unique} $\a<\aleph_{\om+1}$ such that $h=^*f_\a$. Consequently, the
space
$X$ is {\em totally quasi ordered} by
$<^*$, namely the following trichotomy holds: 
\begin{equation}\label{trich}
\forall h,k\in X\;\Bigl[\, h<^* k \;\vee\: k< ^* h\;\vee\;h=^*k \,\Bigr]
\end{equation}

 Claim \ref{suf}  translates  to a  property of $X$:

\begin{claim} \label{convergence}  Suppose that $0<m\le k<\om$ and that
$\lng h_\z:\z<\aleph_m\rng$ is a sequence of elements of $X$ such that
$\lng h_\z\rest (k,\om):\z<\aleph_m\rng$ is increasing in $<$. Denote
$g=\sup\{h_\z:\z<\aleph_m\}$. Then there is some $h\in X$
such that $h=^*g$
\end{claim} 

\begin{proof} For every $\z<\aleph_m$ there is a unique
$\a(\z)<{\aleph_{\om+1}}$ for which $h_\z=^* f_{\a(\z)}$. Since $\lng
h_\z:\z<\aleph_m\rng$ is increasing in $<^*$, the sequence $\lng
\a(\z):\z<\aleph_m\rng$ is strictly increasing. Let
$\d=\sup\{\a(\z):\z<\aleph_m\}$.  By Claim
\ref{suf},   $\cf g(n)=\aleph_m$ for all $n\in (k,\om)$  and $g=^*f_\d$. 

Let $h\in \prod_{n>1}(\aleph_n+1)$ be defined by $h(n)=\aleph_n$ for $n\le
k$ and
$h(n)=g(n)$ for $n>k$. Then $h\in X^R$ and 
$h=^*f_\d$. Thus $h\in X$ and $h=^*g$ as required. 
\end{proof}

\begin{claim} \label{closed} $X$ is a closed subspace of $X^R$. 
\end{claim}

\begin{proof} 
Suppose $t\in \cl X$ and $t\in X^R$. For every $h\in X$ let
$E(h,t):=\{n>1:h(n)=t(n)\}$. 

\begin{claim}\label{fincofin}  If $h\le t$ and $h\in X$ then $E(h,t)$
is either finite or co-finite.
\end{claim}

\begin{proof}[Proof of Claim] Suppose to the contrary that $h\le t$,
$h\in X$ and $|E(h,t)|=|\om - E(h,t)|=\aleph_0$. Let, for $n>1$,
\[
f(n)=
\cases{
0 & if $n\in E(h,t)$\cr
& \cr
h(n)& if $n\in (\om-E(h,t))$\cr}
\]

 Clearly $f<t$. We argue that $X\cap (f,t]$ is empty, contrary to $t\in
\cl X$. Indeed, if $k\in X$ and $k(n)>h(n)$ for all $n\in (w-E(h,t))$
then $k\not <^*h\,\wedge\, k\not=^*h$ because $w-E(h,t)$ is infinite and
so $h<^*k$ by the trichotomy (\ref{trich}).  Since $E(h,t)$ is
infinite and $\{n>1:k(n)\le h(n)\}$ is finite, there is $n\in E(h,t)$
such that $k(n)>h(n)=t(n)$ and therefore  $k\notin (f,t]$.
\end{proof}

We  need a definition:

\begin{definition} ${\bf W}:=\{w\su \om: \forall f<t\;\exists h\in
(f,t]\,
;\Bigl[ E(h,t)=w\Bigr]\}$
\end{definition}

By Claim \ref{fincofin} if $w\in {\bf W}$ then $w$
is finite or $w$ is co-finite. 

\begin{claim} $W\not=\emptyset$
\end{claim}

\begin{proof}[Proof of Claim] Assume that ${\bf W}$ is empty. This is
equivalent, by Claim \ref{fincofin}, to assuming that every finite and
every co-finite
$w\su
\om$ is not in ${\bf W}$. For every finite or co-finite $w\su \om$ fix a
function $f_w<t$ such that $h\in (f_w,t]\cap X\imply E(h,t)\not=w$. Let
$f$ be the supremum of  $f_w$  taken over for all finite and co-finite
$w\su\om$. Since there are countably many $f_w$ and $\cf t(n)>\aleph_0$
for all
$n>1$ it follows that $f<t$. If $h\le t$ is in $X$ and  $w=E(h,t)$ then
$h\notin (f_w,t]$ and hence $h\notin (f,t]$. Thus
$(f,t]\cap X=\emptyset$, contrary to $t\in\cl X$.
\end{proof}

Let us denote $M_m=\{n>1:\cf t(n)=\aleph_m\}$. Likewise, $M_{<
m}=\bigcup _{1<i< m}M_i$. 

\begin{claim} If there is $h\in X$ so that $E(h,t)$ is co-finite then $t\in
X$.
\end{claim}

\begin{proof} Clear.
\end{proof}

\begin{claim} There exists $h\in X$ so that  $E(h,t)$ is co-finite.
\end{claim}

\begin{proof} By Claim \ref{fincofin} is suffices to prove that there is
$h\le t$ in $X$ with  infinite $E(h,t)$.  Let $m$ be the least such that
$M_m$ is infinite. Such $m$ must exist because $\{\cf t(n):n>1\}$ is
bounded by the definition of
$X^R$.   

Fix $w\in {\bf W}$. If $w$ is infinite then we are done. So suppose that
$w$ is finite and let $k=\max\{m,\max w\}$. 

For every $n\in M_m$ fix an increasing sequence $\lng
\g^n_\z:\a<\aleph_m\rng$ with supremum $t(n)$. By induction of
$\z<\aleph_m$ find a sequence $\lng h_\z:\z<\aleph_m\rng$ so that:
\begin{enumerate}
\item $h_\z\le t$ is in $X$ and $E(h_\z,t)=w$
\item $\xi<\z<\aleph_m\imply h_\xi\rest (k,\om)<g h_\z\rest
(k,\om)<t\rest (k,\om)$
\item $h_\z(n)\ge \g_\z^n$ for all $n\in (k,\om)\cap M_m$
\end{enumerate}

At stage $\z$ let $f=\sup\{h_\xi\rest (k,\om):\xi<\z\}$. Since for every
$\xi<\z$ it follows by $E(h_\xi,t)=w$ that $h_\xi\rest (k,\om)<t\rest
(k,\om)$, and since $\cf t(n)\ge \aleph_m$ for all $n\in (k,\om)$, we have
$f<t$. By definition of $w\in {\bf W}$ we can find $h_\z\le t$ in $X$ with
$E(h_\z,t)=w$ such that $h_\z\rest (k,\om)>f\rest (k,\om)$. Without loss
of generality we can choose $h_\z$ so that $h_\z(n)>\g^n_\z$ for all $n\ge
k$ in $M_m$. 

By Claim \ref{convergence} there is some $h\in X$ with $h(n)=^*\sup
\{h_\z(n):\z<\aleph_m\}$. In particular, $h(n)=t(n)$
for all but finitely many  $n\ge k$ in $M_m$. Since $M_m$ is infinite,
$E(h,t)$ is infinite, and we are done. 
\end{proof}
\end{proof}

\begin{claim} \label{normal} $X$ is collectionwise normal
\end{claim}
\begin{proof} Clear from Claim \ref{closed} and Theorem \ref{xrdowker}.
\end{proof}

We  show next that $X$ is not countably paracompact. 

Let $D_n^X=\{f\in X:\exists m\ge n\left[f(m)=\aleph_m\right]\}$ for
$n>1$. It is straightforward that $D_n^X$ is closed and that $\bigcap_n
D_n^X=\emptyset$.

\begin{claim} \label{dowker} If $U_n\su X$ is open,  and $D_n^X\su U_n$
for all $n>1$, then $\bigcap U_n$ is not empty.
\end{claim}

\begin{proof}[Proof of Claim \ref{dowker}]

Suppose that $U_n\supseteq D_n^X$ is open for $n>1$. We need to prove
that $\bigcap_n U_n$ is not empty. 

We shall prove that there is some $f\in \prod_{n>1}\aleph_n$ such that
every $h>f$ in $X$ belongs to this intersection. 

It suffices to show that for each $n>1$  there is some $f_n\in
\prod_{n>1}\aleph_n$ such that $\forall h\in X \,\Bigl[h>f_n\imply f\in
U_n\Bigr]$, because then
$f=\sup\{f_n:1<n<\om\}$ is as required.

Suppose to the contrary that $m>1$ is fixed and for every function $f\in
\prod _{n>1}\aleph_n$ there is some function $h_f>f$ in $X-U_m$. Since
$h_f\notin D_m$, it follows that $h_f(n)<\aleph_n$ for all $n\ge m$.

For a given $f$, let $g_f=\sup\{h_{f'}:f'\in
\prod_{n>1}\aleph_n \wedge (m,\om)\su E(f',f)\}$. Since this supremum is
taken over $\aleph_m$ many functions $h_{f'}$, it follows from the
above that $g_f(n)<\aleph_n$ for all $n> m$. Also. clearly
$g_f(i)=\aleph_i$ for $1<i\le m$.  

Let $\lng M_\z:\z\le\om_1\rng$ be an elementary chain of submodels of
$H(\theta)$ for large enough regular $\theta$ so that:
\begin{itemize} 
\item $\ov f$, $X$ and the functions $f\mapsto h_f$ and $f\mapsto g_f$
belong to $M_0$
\item $M_\z$ has cardinality $\aleph_1$ and $\lng M_\xi:\xi<\z\rng\in
M_{\z+1}$ for all $\zeta<\om_1$. 
\end{itemize}

For every $\z$ let $\chi_\z(n):=\sup (M_\z\cap \aleph_n)$ for all $n>1$.
Since $|M_\z|=\aleph_1$, it follows that $\chi_\z(n)<\aleph_n$ for all
$n$ and hence $\chi_\z\in \prod_{n>1}\aleph_n$. 

Since $\chi_\xi\in M_\z$ for $\xi<\z<\om_1$, by elementarity also
$h_{\chi_\xi}$ and $g_{\chi_\xi}$ belong to $M_\z$ and consequently
$h_{\chi_\xi},g_{\chi_\xi}<\chi_\z$. 

Therefore, if $\xi<\z<\om_1$ then $\chi_\xi <
h_{\chi_\xi}<\chi_\z<h_{\chi_\z}<\chi_{\om_1}$. Thus $\lng
h_{\chi_\z}:\z<\om_1\rng$ is a sequence in $X$, increasing in $<$ with
supremum $\chi_{\om_1}$. By Claim \ref{convergence}, $\chi_{\om_1}\in
X$. 

 Let
$\chi'$ be so that $\chi'(n)=\chi_{\om_1}(n)$ for all $n>m$
 and $\chi'(i)=\aleph_i$ for $1<i\le m$.  So
$\chi'\in D_m^X\su U_m$ and therefore $(f,\chi']\su U_n$ for some
$f<\chi'$, as $U_m$ is open. 

Find some $\z<\om_1$ such that $f\rest (m,\om)<\chi_\z\rest
(m,\om)$. Let 
$f'=f\rest (m+1)\cup \chi_\z\rest (m,\om)$. By the definition of
$g_{\chi_\z}$ we see that
$f'<h_{f'}\le g_{\chi_\z}\le\chi'$ and, of course,
$h_{f'}\notin U_m$.  This contradicts $h_{f'}\in (f,\chi']\su U_m$.
\end{proof}

The space $X$ defined in \ref{X} is normal and not countably
paracompact  by Claim \ref{normal} and Claim \ref{dowker} respectively,
and is therefore  Dowker by \cite{dowker}.
 Since $|X|=\aleph_{\om+1}$ we have proved:

\begin{theorem} There is a ZFC Dowker space of cardinality
$\aleph_{\om+1}$.
\end{theorem}

It is straightforward to verify that the space $X$ constructed above
has weight $\aleph_{\om+1}$ and character $\aleph_\om$.

\begin{problem} Is $\aleph_{\omega+1}$ the {\em first} cardinal in
which one can prove the existence of a Dowker space in ZFC?
\end{problem}

\end{document}